\documentclass[12pt]{article}
\usepackage{amsfonts}

\usepackage{mathrsfs}
\usepackage{amscd}
\usepackage{amsmath,amsfonts,amssymb,amscd}
\usepackage{indentfirst,graphics,epsfig,psfrag}
\input{epsf}

\usepackage{amsmath,amssymb,amsthm, amscd, amsfonts, graphicx,ifpdf}
\usepackage{lineno}

\newtheorem{thm}{Theorem}[section]
\newtheorem{prop}[thm]{Proposition}

\newtheorem{lem}{Lemma}[section]

\def\qed{\nopagebreak\hfill{\rule{4pt}{7pt}}
\medbreak}

\def\pf{\noindent {\it Proof.} }

\title{\bf The minimal size of a graph\\ with generalized connectivity
$\kappa_{3}= 2$\footnote{Supported by NSFC and the Fundamental
Research Funds for the Central Universities. }}

\author{
\small Shasha Li, Xueliang Li, Yongtang Shi\\
\small Center for Combinatorics and LPMC-TJKLC\\
\small Nankai University, Tianjin 300071, China.\\
\small  Email: lss@cfc.nankai.edu.cn, lxl@nankai.edu.cn, shi@nankai.edu.cn\\
}
\date{}
\begin{document}

\maketitle

\begin{abstract}
Let $G$ be a nontrivial connected graph of order $n$ and $k$
an integer with $2\leq k\leq n$. For a set $S$ of $k$ vertices of
$G$, let $\kappa (S)$ denote the maximum number $\ell$ of
edge-disjoint trees $T_1,T_2,\ldots,T_\ell$ in $G$ such that
$V(T_i)\cap V(T_j)=S$ for every pair $i,j$ of distinct integers with
$1\leq i,j\leq \ell$. Chartrand et al. generalized the concept of
connectivity as follows: The $k$-$connectivity$, denoted by
$\kappa_k(G)$, of $G$ is defined by
$\kappa_k(G)=$min$\{\kappa(S)\}$, where the minimum is taken over
all $k$-subsets $S$ of $V(G)$. Thus $\kappa_2(G)=\kappa(G)$, where
$\kappa(G)$ is the connectivity of $G$.

This paper mainly focuses on the minimal number of edges of a graph
$G$ with $\kappa_{3}(G)= 2$. For a graph $G$ of order $v(G)$ and
size $e(G)$ with $\kappa_{3}(G)= 2$, we obtain that $e(G)\geq
\frac{6}{5}v(G)$, and the lower bound
is sharp by showing a class of examples attaining the lower bound.\\[3mm]
\end{abstract}

\section{Introduction}

We follow the terminology and notation of \cite{Bondy} and all
graphs considered here are always simple. As usual, we denote the
numbers of vertices and edges in $G$ by $v(G)$ and $e(G)$, and these
two basic parameters are called the $order$ and $size$ of $G$,
respectively. Let $X$ be a set of vertices of $G$ and $G[X]$ the
subgraph of $G$ whose vertex set is $X$ and whose edge set consists
of all edges of $G$ which have both ends in $X$. A stable set in a
graph is a set of vertices no two of which are adjacent.
The $connectivity$
$\kappa(G)$ of a graph $G$ is defined as the minimum cardinality of
a set $Q$ of vertices of $G$ such that $G-Q$ is disconnected or
trivial. A well-known theorem of Whitney \cite{Whitney} provides an
equivalent definition of the connectivity. For each $2$-subset
$S=\{u,v\}$ of vertices of $G$, let $\kappa(S)$ denote the maximum
number of internally disjoint $uv$-paths in $G$. Then
$\kappa(G)=$min$\{\kappa(S)\}$, where the minimum is taken over all
$2$-subsets $S$ of $V(G)$.

In \cite{Chartrand}, the authors generalized the concept of
connectivity. Let $G$ be a nontrivial connected graph of order $n$
and $k$ an integer with $2\leq k\leq n$. For a set $S$ of $k$
vertices of $G$, let $\kappa (S)$ denote the maximum number $\ell$
of edge-disjoint trees $T_1,T_2,\ldots,T_\ell$ in $G$ such that
$V(T_i)\cap V(T_j)=S$ for every pair $i,j$ of distinct integers with
$1\leq i,j\leq \ell$ (note that the trees are vertex-disjoint in
$G\backslash S$). A collection $\{T_1,T_2,\ldots,T_\ell \}$ of trees
in $G$ with this property is called an {\it internally disjoint set
of trees connecting $S$}. The $k$-$connectivity$, denoted by
$\kappa_k(G)$, of $G$ is then defined by
$\kappa_k(G)=$min$\{\kappa(S)\}$, where the minimum is taken over
all $k$-subsets $S$ of $V(G)$. Thus, $\kappa_2(G)=\kappa(G)$.

In \cite{LLZ}, we focused on the investigation of $\kappa_{3}(G)$
and mainly studied the relationship between the $2$-connectivity and
the $3$-connectivity of a graph. We gave sharp upper and lower
bounds for $\kappa_{3}(G)$ for general graphs $G$, and showed that
if $G$ is a connected planar graph, then $\kappa(G)-1\leq
\kappa_{3}(G)\leq\kappa(G)$. Moreover, we studied the algorithmic
aspects for $\kappa_{3}(G)$ and gave an algorithm to determine
$\kappa_{3}(G)$ for a general graph $G$.

In this paper, we will turn to determining the minimal number of
edges of a graph $G$ with $\kappa_{3}= 2$. For a graph $G$ of order
$v(G)$ and size $e(G)$ with $\kappa_{3}(G)= 2$, we obtain that
$e(G)\geq \frac{6}{5}v(G)$, and the lower bound is sharp by
constructing a class of graphs which attain the lower bound. Note
that for a graph $G$ of order $v(G)$ and size $e(G)$ with
$\kappa(G)= 2$, we only have $e(G)\geq v(G)$, and a cycle of this
order attains the lower bound.

\section{Lower bound}

Before proceeding, we recall a result in \cite{LLZ}, which will be used
frequently in the sequel.

\begin{lem}\label{lem1}
If $G$ is a connected graph with minimum degree $\delta$, then
$\kappa_3(G)\leq \delta$. In particular, if there are two adjacent
vertices of degree $\delta$, then $\kappa_3(G)\leq \delta-1$.
\end{lem}

Now we give the lower bound.

\begin{prop}\label{prop}
Every graph $G$ of order $n$ with $\kappa_3(G)= 2$ has at least
$\frac{6}{5} n$ edges.
\end{prop}
\pf Since $\kappa_3(G)= 2$, by Lemma \ref{lem1}, we know that
$\delta(G)\geq 2$ and any two vertices of degree $2$ are not
adjacent. Denote by $X$ the set of vertices of degree $2$. By Lemma
\ref{lem1}, we have that $X$ is a stable set. Put $Y=V(G)-X$ and
obviously there are $2|X|$ edges joining $X$ to $Y$. Assume that
$m'$ is the number of edges joining two vertices belonging to $Y$.
It is clear that
\begin{align}
e=2|X|+m'.   \label{1}
\end{align}
Since every vertex of $Y$ has degree at least $3$ in $G$, then
$\sum_{v\in Y} d(v)=2|X|+2m' \geq 3|Y|=3(n-|X|)$, namely,
\begin{align}
5|X|+2m' \geq 3n  .  \label{2}
\end{align}
Combining \eqref{1} with \eqref{2}, we have
$\frac{5}{2}e=\frac{5}{2}(2|X|+m')=5|X|+\frac{5}{2}m'\geq 5|X|+2m'
\geq 3n$, namely, $e\geq \frac{6}{5}n$. The proof is complete. \qed

\noindent \textbf{Remark 2.1:} Furthermore, in Proposition
\ref{prop} equality holds if and only if $5|X|+\frac{5}{2}m'=
5|X|+2m'= 3n$, namely, if and only if

\noindent(A) $m'=0$, that is, $Y$ is a stable set and

\noindent(B) the maximum degree $\Delta$ is $3$.\\[2mm]
Moreover, when equality holds, inequality \eqref{2} becomes
$5|X|=3n$, that is, $|X|=\frac{3}{5}n$.

\noindent \textbf{Remark 2.2:} Obviously, for any graph $G$ with
$e(G)=\frac{6}{5}v(G)$, $\kappa_{3}(G)\leq 2$. The next lemma shows
that the number $e(G)=\frac{6}{5}v(G)$ cannot guarantee that
$\kappa_{3}(G)= 2$.
\begin{lem}\label{lem3}
For any connected graph $G$ of order $10$ and size $12$, $\kappa_{3}(G)=1$.
\end{lem}
\pf Note that $e(G)=\frac{6}{5}v(G)$ and so $\kappa_{3}(G)\leq 2$.
Assume, to the contrary, that there is a connected graph $G$ of
order $10$ and size $12$ with $\kappa_{3}(G)=2$. Therefore by Remark
$2.1$, both $X$ and $Y$ are stable sets, $|X|=\frac{3}{5}v(G)=6$ and
$|Y|=4$, where $X$ and $Y$ are the sets of vertices of degrees $2$
and $3$, respectively. Let $X=\{x_{1},\ldots,x_{6}\}$ and
$Y=\{y_{1},\ldots,y_{4}\}$.

\noindent \textbf{Case 1:} For every two vertices $y_{i}$ and $y_{j}$ in $Y$,
there is a vertex in $X$ that is adjacent to both $y_{i}$ and $y_{j}$,
where $1\leq i\neq j\leq 6$.

Note that every vertex in $X$ has degree $2$ and there are exactly six
$2$-subsets of $Y$, namely
$$
\{y_{1},y_{2}\},\{y_{1},y_{3}\},\{y_{1},y_{4}\},\{y_{2},y_{3}\},\{y_{2},y_{4}\},\{y_{3},y_{4}\}.
$$
Thus we may assume that $G$ is isomorphic to Figure $1$. Then observe
that it is impossible to find two internally-disjoint trees connecting the
vertices $x_{1}$, $x_{2}$ and $x_{4}$, contrary to our assumption.

\begin{center}
\begin{picture}(120,90)
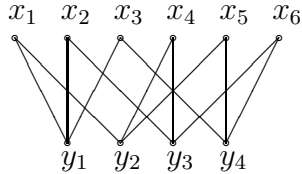
\put(-60,0){Figure $1$: The graph for Case $1$ of Lemma 2.2}

\put(82.5,22){$y_4$}
\put(22.5,22){$y_1$}
\put(42.5,22){$y_2$}
\put(62.5,22){$y_3$}

\put(2.5,77){$x_1$}
\put(22.5,77){$x_2$}
\put(42.5,77){$x_3$}
\put(62.5,77){$x_4$}
\put(82.5,77){$x_5$}
\put(102.5,77){$x_6$}

\put(85,30){\circle{2}}
\put(25,30){\circle{2}}
\put(45,30){\circle{2}}
\put(65,30){\circle{2}}

\put(85,70){\circle{2}}
\put(25,70){\circle{2}}
\put(45,70){\circle{2}}
\put(65,70){\circle{2}}
\put(5,70){\circle{2}}
\put(105,70){\circle{2}}

\put(5,70){\line(1,-2){20}}
\put(5,70){\line(1,-1){40}}
\put(25,70){\line(0,-1){40}}
\put(25,70){\line(1,-1){40}}
\put(45,70){\line(-1,-2){20}}
\put(45,70){\line(1,-1){40}}
\put(65,70){\line(-1,-2){20}}
\put(65,70){\line(0,-1){40}}
\put(85,70){\line(-1,-1){40}}
\put(85,70){\line(0,-1){40}}
\put(105,70){\line(-1,-1){40}}
\put(105,70){\line(-1,-2){20}}

\end{picture}
\end{center}

\noindent \textbf{Case 2:} For some two vertices $y_{i}$ and $y_{j}$ in $Y$,
at least two vertices in $X$ are adjacent to both $y_{i}$ and $y_{j}$,
where $1\leq i\neq j\leq 6$. Since
$G$ is connected, we can get that only two vertices in $X$ are adjacent to both
$y_{i}$ and $y_{j}$. Then we may assume that $G$ is isomorphic to Figure $2$.
Now consider the three vertices $x_{1}$, $x_{3}$ and $x_{5}$ and we can get
$\kappa_{3}(G)=1$, contrary to our assumption.

\begin{center}
\begin{picture}(120,90)
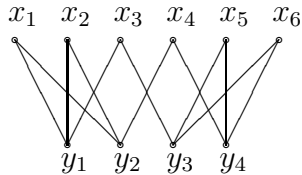
\put(-60,0){Figure $2$: The graph for Case $2$ of Lemma 2.2}

\put(82.5,22){$y_4$}
\put(22.5,22){$y_1$}
\put(42.5,22){$y_2$}
\put(62.5,22){$y_3$}

\put(2.5,77){$x_1$}
\put(22.5,77){$x_2$}
\put(42.5,77){$x_3$}
\put(62.5,77){$x_4$}
\put(82.5,77){$x_5$}
\put(102.5,77){$x_6$}

\put(85,30){\circle{2}}
\put(25,30){\circle{2}}
\put(45,30){\circle{2}}
\put(65,30){\circle{2}}

\put(85,70){\circle{2}}
\put(25,70){\circle{2}}
\put(45,70){\circle{2}}
\put(65,70){\circle{2}}
\put(5,70){\circle{2}}
\put(105,70){\circle{2}}

\put(5,70){\line(1,-2){20}}
\put(5,70){\line(1,-1){40}}
\put(25,70){\line(0,-1){40}}
\put(25,70){\line(1,-2){20}}
\put(45,70){\line(-1,-2){20}}
\put(45,70){\line(1,-2){20}}
\put(65,70){\line(-1,-2){20}}
\put(65,70){\line(1,-2){20}}
\put(85,70){\line(-1,-2){20}}
\put(85,70){\line(0,-1){40}}
\put(105,70){\line(-1,-1){40}}
\put(105,70){\line(-1,-2){20}}

\end{picture}
\end{center}

The proof is complete. \qed

Next we will show that the lower bound given in Proposition
\ref{prop} is essentially best possible. For this, we construct a
class of graphs attaining the lower bound.

Before proceeding, we want to give some notions. For any two
integers $a$ and $k\geq 1$, denote by $[a]_{k}$ an integer such that
$1\leq [a]_{k}\leq k$ and $a\equiv [a]_{k}$ $(mod\ k)$. For a cycle
$C=x_{1}x_{2}x_{3}\ldots x_{k-1}x_{k}x_{1}$, we denote three special
segments of $C$ by
$x_{a}Cx_{b}=x_{a}x_{[a+1]_{k}}x_{[a+2]_{k}}$ $\ldots
x_{[b-1]_{k}}x_{b}$,
$\hat{x}_{a}Cx_{b}=x_{[a+1]_{k}}x_{[a+2]_{k}}\ldots
x_{[b-1]_{k}}x_{b}$ and
$\hat{x}_{a}C\hat{x}_{b}=x_{[a+1]_{k}}x_{[a+2]_{k}}\ldots
x_{[b-1]_{k}}$, where $1\leq a, b\leq k$. Denote by
$|C|$ and $|P|$ the lengths of a cycle $C$ and a path $P$,
respectively.

\begin{lem}\label{lem2}
For a positive integer $k\neq 2$, let $C=x_{1}y_{1}x_{2}y_{2}\ldots
x_{2k}y_{2k}x_{1}$ be a cycle of length $4k$. Add $k$ new vertices
$z_{1},z_{2},\ldots,z_{k}$ to $C$, and join $z_{i}$ to $x_{i}$ and
$x_{i+k}$, for $1\leq i\leq k$. The resulting graph is denoted by
$H$. Then, the $3$-connectivity of $H$ is $2$, namely,
$\kappa_3(H)=2$.
\end{lem}
\pf Since $\delta(H)=2$, by Lemma \ref{lem1} we can get
$\kappa_3(H)\leq 2$. So the task is to show $\kappa_3(H)\geq 2$. By
the definition of the generalized connectivity, it suffices to prove
that $\kappa(S)\geq 2$, for every $3$-subset $S$ of $V(H)$.

Firstly, partition $V(H)$ into three types:
$V_{1}=\{x_{1},x_{2},\ldots,x_{2k}\}$,
$V_{2}=\{z_{1},z_{2},\ldots,z_{k}\}$ and
$V_{3}=\{y_{1},y_{2},\ldots,y_{2k}\}$. We proceed by considering all
cases of $S$.

\textbf{Case 1:} $S=\{x_{a},x_{b},x_{c}\}$, where $1\leq a< b< c\leq 2k$.

The three vertices divide the cycle $C$ into three segments, at
least one of which has length at most $|C|/3$. Without loss of
generality, we may assume that $|x_{a}Cx_{b}|\leq |C|/3$, namely,
$|x_{b}Cx_{a}|\geq 2|C|/3$. Let $b'=[b+k]_{2k}$. Note that
$|x_{b}Cx_{b'}|=|C|/2$, and so $x_{b'}\in
V(\hat{x}_{b}C\hat{x}_{a})$.

\textbf{Subcase 1.1:} $x_{b'}\in V(x_{c}C\hat{x}_{a})$. In this
case, $T_{1}=x_{a} Cx_{b}Cx_{c}$ and $T_{2}=x_{c}Cx_{b'}Cx_{a}\cup
x_{b'}z_{[b]_{k}}x_{b}$ are two internally disjoint trees connecting
$S$.

\textbf{Subcase 1.2:} $x_{b'}\in V(\hat{x}_{b}C\hat{x}_{c})$. Let
$a'=[a+k]_{2k}$. We can get $x_{a'}\in V(\hat{x}_{b}C\hat{x}_{b'})$,
since $1\leq |x_{a}Cx_{b}|\leq |C|/3$, $|x_{a}Cx_{a'}|=|C|/2$ and
$|x_{b}Cx_{b'}|=|C|/2$. Therefore, $x_{a'}\in
V(\hat{x}_{b}C\hat{x}_{c})$, and then $T_{1}=x_{c}Cx_{a}Cx_{b}$ and
$T_{2}=x_{b}Cx_{a'}Cx_{c}\cup x_{a'}z_{[a]_{k}}x_{a}$ are two
internally disjoint trees connecting $S$.

\textbf{Case 2:} $S=\{z_{a},z_{b},z_{c}\}$, where $1\leq a< b< c\leq k$.

Since $1\leq a< b< c\leq k <a+k< b+k <c+k\leq 2k$, $x_{a}Cx_{b}Cx_{c}$
and $x_{a+k}Cx_{b+k}Cx_{c+k}$ are two disjoint segments of $C$. It is easy to find two
internally disjoint trees connecting $S$:
$T_{1}=z_{a}x_{a}Cx_{b}Cx_{c}z_{c}\cup x_{b}z_{b}$ and $T_{2}=z_{a}x_{a+k}Cx_{b+k}Cx_{c+k}z_{c}\cup x_{b+k}z_{b}$.

\textbf{Case 3:} $S=\{x_{a},x_{b},z_{c}\}$, where $1\leq a< b\leq 2k$ and $1\leq c\leq k$.

Observe that the two neighbors $x_{c}$ and $x_{c+k}$ of $z_{k}$
divide the cycle into two segments $x_{c}Cx_{c+k}$ and
$x_{c+k}Cx_{c}$.

\textbf{Subcase 3.1:} $x_{a}$ and $x_{b}$ lie in distinct segments.
Without loss of generality, we may assume that $x_{a}\in
V(x_{c}Cx_{c+k})$ and $x_{b}\in V(x_{c+k}Cx_{c})$. Now
$T_{1}=x_{a}Cx_{c+k}Cx_{b}\cup x_{c+k}z_{c}$ and
$T_{2}=x_{b}Cx_{c}Cx_{a}\cup x_{c}z_{c}$ are two trees we want. Note
that the subcase contains the situation that either $x_{c}$ or
$x_{c+k}$ is exactly $x_{a}$ or $x_{b}$.

\textbf{Subcase 3.2:} $x_{a}$ and $x_{b}$ lie in the same segment.
Without loss of generality, suppose that $x_{a},x_{b}\in
V(\hat{x}_{c}C\hat{x}_{c+k})$. Let $b'=[b+k]_{2k}$. Since
$|x_{c}Cx_{c+k}|=|C|/2$, $|x_{b}Cx_{b'}|=|C|/2$ and $x_{b}\in
V(\hat{x}_{c}C\hat{x}_{c+k})$, we have $x_{b'}\in
V(\hat{x}_{c+k}C\hat{x}_{c})$ and $T_{1}=x_{a}Cx_{b}Cx_{c+k}z_{c}$
and $T_{2}=x_{b}z_{[b]_{k}}x_{b'}Cx_{c}Cx_{a}\cup x_{c}z_{c}$ are
two internally disjoint trees connecting $S$.

\textbf{Case 4:} $S=\{x_{a},z_{b},z_{c}\}$, where $1\leq a\leq 2k$ and $1\leq b< c\leq k$.

Since $1\leq b< c \leq k <b+k <c+k \leq 2k$, the two neighbors
$x_{b},x_{b+k}$ of $z_{b}$, together with two neighbors
$x_{c},x_{c+k}$ of $z_{c}$ divide the cycle into four segments
$x_{b}Cx_{c}$, $x_{c}Cx_{b+k}$, $x_{b+k}Cx_{c+k}$ and
$x_{c+k}Cx_{b}$. Actually, it is easy to see that no matter which
segment $x_{a}$ lies in, the situations are equivalent. Therefore, without
loss of generality, we may assume that $x_{a}\in V(x_{b}Cx_{c})$. We
have $T_{1}=x_{a}Cx_{c}Cx_{b+k}z_{b}\cup x_{c}z_{c}$ and
$T_{2}=z_{c}x_{c+k}Cx_{b}Cx_{a}\cup x_{b}z_{b}$ are two internally
disjoint trees connecting $S$. Note that this case includes the
situation that $x_{a}$ is exactly $x_{b}$ or $x_{c}$.

Next we consider the cases in which $S$ contains the vertices in
$V_3$.

\textbf{Case 5:} $S=\{y_{a},y_{b},y_{c}\}$, where $1\leq a< b< c\leq 2k$.

Clearly, in this case, $k$ is a positive integer at least $3$. Among
the three segments $y_{a}Cy_{b}$, $y_{b}Cy_{c}$ and $y_{c}Cy_{a}$ of
$C$, at least one of them has length not more than $|C|/3$. We may
assume that $|y_{a}Cy_{b}|\leq |C|/3=4k/3$. Moreover, observe that
$x_{a+1}$ lies between $y_{a}$ and $y_{b}$. We have $y_{b}\in
V(\hat{x}_{a+1}C\hat{x}_{[a+1+k]_{2k}})$, since
$|x_{a+1}Cy_{b}|<|y_{a}Cy_{b}| \leq 4k/3$ and
$|x_{a+1}Cx_{[a+1+k]_{2k}}|=|C|/2=2k$.

\textbf{Subcase 5.1:} $y_{c}\in
V(\hat{y}_{b}C\hat{x}_{[a+1+k]_{2k}})$. There is at least one vertex
$x_{b+1}$ between $y_{b}$ and $y_{c}$. Since $x_{b+1}\in
V(\hat{x}_{a+1}C\hat{x}_{[a+1+k]_{2k}})$, it is clear that
$x_{[b+1+k]_{2k}}\in V(\hat{x}_{[a+1+k]_{2k}}C\hat{x}_{a+1})$,
namely, $x_{[b+1+k]_{2k}}\in V(\hat{x}_{[a+1+k]_{2k}}C\hat{y}_{a})$.
We can find two internally disjoint trees connecting $S$:
$T_{1}=y_{a}x_{a+1}Cy_{b}\cup y_{c}Cx_{[a+1+k]_{2k}}\cup
x_{a+1}z_{[a+1]_{k}}x_{[a+1+k]_{2k}}$ and
$T_{2}=y_{b}x_{b+1}Cy_{c}\cup
x_{b+1}z_{[b+1]_{k}}x_{[b+1+k]_{2k}}Cy_{a}$.

\textbf{Subcase 5.2:} $y_{c}\in
V(\hat{x}_{[a+1+k]_{2k}}C\hat{y}_{a})$. There is at least one vertex
$x_{a}$ between $y_{c}$ and $y_{a}$. Obviously, $x_{[a+k]_{2k}}\in
V(\hat{x}_{a+1}C\hat{x}_{[a+1+k]_{2k}})$. Moreover,
$x_{a}Cy_{b}=|y_{a}Cy_{b}|+1\leq |C|/3+1=4k/3+1$ and
$x_{a}Cx_{[a+k]_{2k}} =|C|/2=2k$, where $k\geq 3$. So $y_{b}\in
V(\hat{x}_{a}C\hat{x}_{[a+k]_{2k}})$. Now
$T_{1}=y_{a}x_{a+1}Cy_{b}\cup
x_{a+1}z_{[a+1]_{k}}x_{[a+1+k]_{2k}}Cy_{c}$ and
$T_{2}=y_{b}Cx_{[a+k]_{2k}}z_{[a]_{k}}x_{a}\cup y_{c}Cx_{a}y_{a}$
are two internally disjoint trees connecting $S$.

\textbf{Case 6:} $S=\{y_{a},y_{b},x_{c}\}$, where $1\leq a< b\leq 2k$ and $1\leq c\leq 2k$.

Notice that $y_{a}$ and $y_{b}$ divide $C$ into two segments
$y_{a}Cy_{b}$ and $y_{b}Cy_{a}$. Let $c'=[c+k]_{2k}$, and then two
subcases arise.

\textbf{Subcase 6.1:} $x_{c}$ and $x_{c'}$ lie in distinct segments.
We may assume that $x_{c}\in V(y_{a}Cy_{b})$ and $x_{c'}\in
V(y_{b}Cy_{a})$. Thus, $T_{1}=y_{a}Cx_{c}Cy_{b}$ and
$T_{2}=y_{b}Cx_{c'}Cy_{a}\cup x_{c}z_{[c]_{k}}x_{c'}$ are exactly
two trees we want.

\textbf{Subcase 6.2:} $x_{c}$ and $x_{c'}$ lie in the same segment.
Without loss of generality, we may assume that $x_{c},x_{c'}\in
V(y_{b}Cy_{a})$ and they occur in cyclic order $y_{a},
y_{b},x_{c},x_{c'}$ on $C$. The segment $y_{a}Cy_{b}$ must contain a
vertex $x_{a+1}$ in $V_{1}$. Since $x_{a+1}\in
V(\hat{x}_{c'}C\hat{x}_{c})$, $x_{[a+1+k]_{2k}}\in
V(\hat{x}_{c}C\hat{x}_{c'})$. So we can find two internally disjoint
trees connecting $S$: $T_{1}=y_{a}x_{a+1}Cy_{b}\cup
x_{a+1}z_{[a+1]_{k}}x_{[a+1+k]_{2k}}\cup x_{c}Cx_{[a+1+k]_{2k}}$ and
$T_{2}=y_{b}Cx_{c}z_{[c]_{k}}x_{c'}Cy_{a}$.

\textbf{Case 7:} $S=\{y_{a},y_{b},z_{c}\}$, where $1\leq a< b\leq 2k$ and $1\leq c\leq k$.

If $k=1$, then $C=x_{1}y_{1}x_{2}y_{2}x_{1}$ and $H=C\cup
x_{1}z_{1}x_{2}$. So $y_{a},y_{b}$ and $z_{c}$ are exactly
$y_{1},y_{2}$ and $z_{1}$, respectively. Now
$T_{1}=y_{2}x_{1}y_{1}\cup x_{1}z_{1}$ and
$T_{2}=y_{1}x_{2}y_{2}\cup x_{2}z_{1}$ are two internally disjoint
trees connecting $S$.

Otherwise, $k\geq 3$, since $k\neq 2$. We know that $y_{a},y_{b}$
divide $C$ into two segments $y_{a}Cy_{b},y_{b}Cy_{a}$, and $z_{c}$
has two neighbors $x_{c}$ and $x_{c+k}$.

\textbf{Subcase 7.1:} $x_{c}$ and $x_{c+k}$ lie in distinct segments. Suppose that
$x_{c}\in V(y_{a}Cy_{b})$ and $x_{c+k}\in V(y_{b}Cy_{a})$. Clearly $T_{1}=y_{a}Cx_{c}Cy_{b}
\cup x_{c}z_{c}$ and $T_{2}=y_{b}Cx_{c+k}Cy_{a}\cup x_{c+k}z_{c}$
are two internally disjoint trees connecting $S$.

\textbf{Subcase 7.2:} $x_{c}$ and $x_{c+k}$ lie in the same segment.
Without loss of generality, we may assume that $x_{c},x_{c+k}\in
V(y_{b}Cy_{a})$ and they occur in cyclic order $y_{a},
y_{b},x_{c},x_{c+k}$ on $C$.

\textbf{Subsubcase 7.2.1:} Between $y_{a}$ and $y_{b}$, there are at
least two vertices in $V_{1}$. Clearly $x_{a+1}\neq x_{b}$, and
$y_{a},x_{a+1},x_{b},y_{b},x_{c},x_{[a+1+k]_{2k}}, x_{[b+k]_{2k}}$
and $x_{c+k}$ are the cyclic order in which they occur on $C$. So
we can find two internally disjoint trees connecting $S$:
$T_{1}=y_{a}x_{a+1}z_{[a+1]_{k}}x_{[a+1+k]_{2k}}\cup
y_{b}Cx_{c}Cx_{[a+1+k]_{2k}}\cup x_{c}z_{c}$ and
$T_{2}=y_{b}x_{b}z_{[b]_{k}}x_{[b+k]_{2k}}Cx_{c+k}Cy_{a}\cup
x_{c+k}z_{c}$.

\textbf{Subsubcase 7.2.2:} Between $y_{a}$ and $y_{b}$, there is
only one vertex in $V_{1}$, i.e, $x_{a+1}= x_{b}$. Let
$b'=[b+k]_{2k}$ and clearly $x_{b'}\in V(\hat{x}_cC\hat{x}_{c+k})$.
Since $k\geq 3$, $V(\hat{x}_cC\hat{x}_{c+k})$ contains at least two
vertices $x_{c+1},x_{c+k-1}$ in $V_{1}$. If $x_{c+1}\neq x_{b'}$,
then $x_{[c+1+k]_{2k}}=x_{[c+k+1]_{2k}}\neq x_{b}\in
V(\hat{x}_{c+k})C\hat{y}_{a}$. So $T_{1}=y_{a}x_{b}y_{b}\cup
x_{b}z_{[b]_{k}}x_{b'}Cx_{c+k}z_{c}$ and
$T_{2}=y_{b}Cx_{c}y_{c}x_{c+1}z_{[c+1]_{k}}x_{[c+k+1]_{2k}}Cy_{a}\cup x_{c}z_{c}$
are two internally disjoint trees connecting $S$. Otherwise,
$x_{c+k-1}\neq x_{b'}$, i.e, $x_{[c-1]_{2k}}\neq x_{b}$. We have
$x_{[c-1]_{2k}}\in V(\hat{y}_{b}C\hat{x}_{c})$. So
$T_{1}=y_{a}x_{b}y_{b}\cup x_{b}z_{[b]_{k}}x_{b'}\cup
z_{c}x_{c}Cx_{b'}$ and
$T_{2}=y_{b}Cx_{[c-1]_{2k}}z_{[c-1]_{k}}x_{c+k-1}y_{c+k-1}x_{c+k}Cy_{a}\cup
x_{c+k}z_{c}$ are two internally disjoint trees connecting $S$.

\textbf{Case 8:} $S=\{y_{a},x_{b},x_{c}\}$, where $1\leq a\leq 2k$ and $1\leq b< c\leq 2k$.

Let $b'=[b+k]_{2k}$ and $c'=[c+k]_{2k}$. If $b'=c$, i.e.,
$c=[b+k]_{2k}$, then without loss of generality, we may assume that
$y_{a}\in V(x_{b}Cx_{c})$. We have $T_{1}=y_{a}Cx_{c}z_{[c]_{k}}x_{b}$
and $T_{2}=x_{c}Cx_{b}Cy_{a}$ are two internally disjoint trees
connecting $S$. Otherwise, $b'\neq c$. Without loss of generality,
suppose $x_{b},x_{c},x_{b'}$ and $x_{c'}$ are the cyclic order in
which they occur on $C$, and then they divide $C$ into four
segments $x_{b}Cx_{c},x_{c}Cx_{b'},x_{b'}Cx_{c'}$ and
$x_{c'}Cx_{b}$.

\textbf{Subcase 8.1:} $y_{a}\in V(x_{b}Cx_{c})$. We can find two
internally disjoint trees connecting $S$: $T_{1}=x_{b}Cy_{a}\cup
x_{c}Cx_{b'}z_{[b]_{k}}x_{b}$ and $T_{2}
=y_{a}Cx_{c}z_{[c]_{k}}x_{c'}Cx_{b}$.

\textbf{Subcase 8.2:} $y_{a}\in V(x_{c}Cx_{b'})$ or $y_{a}\in
V(x_{c'}Cx_{b})$. It is easy to see that the two situations are
actually equivalent. So we only consider the former. We can find two
internally disjoint trees connecting $S$: $T_{1}=x_{b}Cx_{c}Cy_{a}$
and $T_{2}=y_{a}Cx_{b'}Cx_{c'}z_{[c]_{k}}x_{c} \cup
x_{b'}z_{[b]_{k}}x_{b}$.

\textbf{Subcase 8.3:} $y_{a}\in V(x_{b'}Cx_{c'})$. We can find two
internally disjoint trees connecting $S$: $T_{1}=x_{b}Cx_{c}\cup
x_{b}z_{[b]_{k}}x_{b'}Cy_{a}$ and $T_{2}=y_{a}Cx_{c'}Cx_{b}\cup
x_{c'}z_{[c]_{k}}x_{c}$.

\textbf{Case 9:} $S=\{y_{a},z_{b},z_{c}\}$, where $1\leq a\leq 2k$ and $1\leq b< c\leq k$.

Observe that $x_{b},x_{c},x_{b+k}$ and $x_{c+k}$ divide the cycle
into four segments $x_{b}Cx_{c},x_{c}Cx_{b+k}$,
$x_{b+k}Cx_{c+k}$ and $x_{c+k}Cx_{b}$. Actually, no matter which
segment $y_{a}$ lies in, the situations are equivalent. So without
loss of generality, we may assume that $y_{a}\in V(x_{b}Cx_{c})$.
Now $T_{1}=y_{a}Cx_{c}Cx_{b+k}z_{b}\cup x_{c}z_{c}$ and
$T_{2}=z_{c}x_{c+k}Cx_{b}Cy_{a}\cup x_{b}z_{b}$ are two internally
disjoint trees connecting $S$.

\textbf{Case 10:} $S=\{y_{a},x_{b},z_{c}\}$, where $1\leq a\leq 2k$, $1\leq b\leq 2k$ and $1\leq c\leq k$.

\textbf{Subcase 10.1:} $b=c$ or $b=c+k$. Without loss of generality,
we may assume that $b=c$ and $y_{a}\in V(x_{c+k}Cx_{b})$. Therefore,
$T_{1}=y_{a}Cx_{b}z_{c}$ and $T_{2}=x_{b}Cx_{c+k}Cy_{a}\cup
x_{c+k}z_{c}$ are two internally disjoint trees connecting $S$.

\textbf{Subcase 10.2:} $b\neq c$ and $b\neq c+k$.
Let $b'=[b+k]_{2k}$. We may assume that $x_{b},x_{c},x_{b'}$ and $x_{c+k}$ are the
cyclic order in which they occur on $C$. Moreover, they divide $C$
into four segments $x_{b}Cx_{c}$, $x_{c}Cx_{b'}$, $x_{b'}Cx_{c+k}$ and $x_{c+k}Cx_{b}$.

If $y_{a}\in V(x_{b}Cx_{c})$, then $T_{1}=y_{a}Cx_{c}Cx_{b'}z_{[b]_{k}}x_{b}\cup x_{c}z_{c}$
and $T_{2}=z_{c}x_{c+k}Cx_{b}Cy_{a}$
are two internally disjoint trees connecting $S$.

If $y_{a}\in V(x_{c}Cx_{b'}Cx_{c+k})$, then $T_{1}=x_{b}Cx_{c}Cy_{a}\cup x_{c}z_{c}$
and $T_{2}=y_{a}Cx_{c+k}Cx_{b}\cup x_{c+k}z_{c}$
are two internally disjoint trees connecting $S$.

If $y_{a}\in V(x_{c+k}Cx_{b})$, then $T_{1}=y_{a}Cx_{b}Cx_{c}z_{c}$
and $T_{2}=x_{b}z_{[b]_{k}}x_{b'}Cx_{c+k}Cy_{a}\cup x_{c+k}z_{c}$
are two internally disjoint trees connecting $S$.

The proof is complete. \qed

\noindent \textbf{Remark 2.3:} Clearly the order $v(H)$ of the graph
$H$ is $5k$ and the size $e(H)$ is $4k+2k=6k$, where $k\neq 2$ is a
positive integer. Therefore $e(H)=\frac{6}{5}v(H)$, and by Lemma
\ref{lem2}, we know that $\kappa_{3}(H)=2$. It follows that $H$
attains the lower bound of Proposition \ref{prop}.

\noindent \textbf{Remark 2.4:} If $k=2$,  then $H$ is a connected
graph of order $10$ and size $12$. By Lemma \ref{lem3}, we can get
$\kappa_{3}(H)=1$. This is the reason why we add the condition
$k\neq 2$ to Lemma \ref{lem2}. Moreover, no graphs of order $10$ can
attain the lower bound.

Now, we can obtain our main result.
\begin{thm}
If $G$ is a graph of order $n$ with $\kappa_{3}(G)= 2$, then
$e(G)\geq \frac{6}{5}n$ and the lower bound is sharp.\qed
\end{thm}

\noindent{Acknowledgement:} The authors would like to thank the
referees for comments and suggestions, which helped to improve the
presentation of the paper.


\begin{thebibliography}{99}

\bibitem{Bondy} J.A. Bondy and U.S.R. Murty, Graph Theory, GTM 244,
 Springer, 2008.

\bibitem{Chartrand} G. Chartrand, F. Okamoto, P. Zhang, Rainbow
trees in graphs and generalized connectivity, Networks 55(4)(2010),
360--367 .

\bibitem{LLZ} S. Li, X. Li, W. Zhou, Sharp bounds for the generalized
connectivity $\kappa_3(G)$, Discrete Math. 310(2010), 2147--2163.

\bibitem{Whitney} H. Whitney, Congruent graphs and the connectivity of graphs
and the connectivity of graphs, Amer. J. Math. 54(1932), 150--168.

\end{thebibliography}
\end{document}